\documentclass{amsart}
\usepackage{amsmath}
\usepackage{amsfonts}
\usepackage{amssymb}
\usepackage{graphicx}
\usepackage{float, verbatim}
\usepackage{enumerate}
\usepackage{color}

\newcommand{\Haus}{\dim_{\mathrm{H}}}
\newcommand{\Assouad}{\dim_\mathrm{A}}
\newcommand{\boxd}{\dim_\mathrm{B}}

\newtheorem{thm}{Theorem}[section]

\newtheorem{ques}[thm]{Question}

\newtheorem{defn}[thm]{Definition}

\newtheorem{rem}[thm]{Remark}
\newtheorem*{thm*}{Theorem}
\newtheorem*{conj*}{Conjecture}
\newtheorem*{lma*}{Lemma}

\begin{document}

\title{Cube packings in Euclidean spaces}

\author{Han Yu}
\address{Han Yu\\
School of Mathematics \& Statistics\\University of St Andrews\\ St Andrews\\ KY16 9SS\\ UK \\ }
\curraddr{}
\email{hy25@st-andrews.ac.uk}
\thanks{}



\date{}

\dedicatory{}

\begin{abstract}
In this paper we study some cube packing problems. In particular we are interested in compact subsets of $\mathbb{R}^n,n\geq 2$, which contain boundaries of cubes with all side lengths in $(0,1)$. We show here that such sets must have lower box dimension at least $n-0.5$ and we will also provide sharp examples. We also show here that such sets must be large in general in a precise sense which is also introduced in this paper.
\end{abstract}

\maketitle

\section{Cube packings}
The central topic of this paper is cube packing in $\mathbb{R}^n, n\geq 2$. Consider the unit cube $K=[-0.5,0.5]^n$, we denote $K_{n-1}=\partial K$ to be its boundary which is a union of $2n$ faces. In general, we use $K_{i},0\leq i\leq n$ for the $i$-skeleton of the unit cube. In this paper we shall only study $K_{n-1}$. In particular, our cubes are aligned with the coordinate axes.

The problem we are interested in is how to pack cubes of different sizes economically in the sense of box dimension. We need to formally define the term 'packing'.

\begin{defn}[Cube packing: different sizes]\label{defsi}
	A size cube packing is defined as a map:
	\[
	f_p:(0,1)\to [0,1]^n
	\]
	
	The image of packing $f_p$ is 
	\[
	T(f_p)=\bigcup_{t\in (0,1)} tK_{n-1}+f_p(t)\subset\mathbb{R}^n.
	\]
	
	We call the packing measurable if $f_p$ is Borel measurable.
\end{defn} 

Intuitively speaking, $f_p(t)$ is the centre of a cube with side length $t$ and $T(f_p)$ is the union of the $n-1$ skeletons of all such cubes. Similarly we can also consider the following notion. 

\begin{defn}[Cube packing: different centres]\label{defpo}
	A centre cube packing is defined as a map:
	\[
	g_p:[0,1]^n\to (0,1)
	\]
	
	The image of packing $g_p$ is 
	\[
	T(g_p)=\bigcup_{x\in [0,1]^n} g_p(x)K_{n-1}+x\subset\mathbb{R}^n.
	\]
	
	We call the packing measurable if $g_p$ is Borel measurable.
\end{defn} 

Such packing problems can be dated back to Bourgain \cite{Bo} and Marstrand \cite{Mas} and later by Wolff \cite{W1} and Schlag \cite{Sch}. They considered sphere packings, namely, $K_{n-1}$ is replaced by $S^{n-1}$ in Definition \ref{defpo}. Now it is known that if $A\subset\mathbb{R}^2$ has Hausdorff dimension larger than $1$, then any set which contains circles with centre set precisely $A$ must have positive measure.

For packing spheres with different sizes, it is known by Kolasa and Wolff \cite{W2} and Wolff \cite{W1} that for any $n\geq 2$, any compact subset of $\mathbb{R}^n$ which contains sphere of all radii in $[1/2,1]$ has full Hausdorff dimension. It turns out that the most difficult case is when $n=2$. Kolasa and Wolff developed a method involving cone constructions and cell decompositions. This method relies heavily on the non-vanishing curvature of circles. 

For cube packing with different centres, such problems were considered by Keleti, Nagy, Shmerkin \cite{SK}, Thornton \cite{Th} and  Chang, Cs\"{o}rnyei, H\'{e}ra, Keleti \cite{CCHK}. In particular in $\mathbb{R}^2$ they showed that a set which contains cubes with all centres in $[0,1]\times [0,1]$ has Hausdorff dimension at least $1$ and lower box dimension at least $7/4$. Their result is sharp. 

From now on we shall focus on cube packings with different sizes (Definition \ref{defsi}).

\begin{thm}\label{Th1}
	For any $n\geq 2$, consider an arbitrary size cube packing $f_p$. Denote $T(f_p)$ as $G$ then we have the following results:
	\[
	\Assouad G\geq\underline{\boxd} G\geq n-0.5.
	\]
	\[
	\Haus G\geq n-1.
	\]
\end{thm}

The above result is sharp and we shall provide example to illustrate this.

\begin{thm}\label{Th1.1}
	For any $n\geq 2$ and $\epsilon>0$, there exist size cubes packings with images $G_1,G_2,G_3$ such that
	\[
	\Haus G_1=n-1.
	\]
	
	\[
	\boxd G_2=n-0.5.
	\]
	
	\[
	\Assouad G_3\leq n-0.5+\epsilon.
	\]
\end{thm}

By knowing the results of sphere packings( They are usually full dimension or positive Lebesgue measure) we do intuitively think images of size cube packings must be in some sense large. Towards this direction we consider the original motivation of circle maximal problems. Given a wave equation, the solution can be presented as 'wave-fronts', if the source set is given, then how large is the 'wave-front' set in relation with time? Then we can consider similar problem with 'cube-front' and the precise formulation is included in the following theorems. 

\begin{thm}\label{Th2}
	For any $n\geq 2$, consider a measurable size cube packing $f_p$. Let $G=T(f_p)$ and for each $r\in (0,1)$ we construct the following set
	\[
	G_r=\bigcup_{t\in (0,1)} rt K_{n-1}+f_p(t).
	\]
	Then for Lebesgue almost all $r\in (0,1)$, $\Haus G_r=n$. 
\end{thm}
\begin{rem}
	Intuitively speaking, $G_r$ is a set obtained by shrinking each cube in $G$ by ratio $r$ and keeping the centre.
\end{rem}

\section{Notations and preliminaries}
Here we list some notions of dimensions that will be used in this paper. We refer \cite[chapter 4,5]{Ma}, \cite[chapter 2,3]{Fa} for more details and properties of these dimensions.

We shall use $N_r(F)$ for the minimal covering number of a set $F$ in $\mathbb{R}^n$ with cubes of side length $r>0$. 

\textbf{Hausdorff dimension:}
For any $s\in\mathbb{R}^+$, for any $\delta>0$ define the following quantity:
\[
\mathcal{H}^s_\delta(F)=\inf\left\{\sum_{i=1}^{\infty}(\mathrm{diam} (U_i))^s: \bigcup_i U_i\supset F, \mathrm{diam}(U_i)<\delta\right\}.
\]

Then the $s$-Hausdorff measure of $F$ is:
\[
\mathcal{H}^s(F)=\lim_{\delta\to 0} \mathcal{H}^s_{\delta}(F).
\]

The Hausdorff dimension of $F$ is:
\[
\Haus F=\inf\{s\geq 0:\mathcal{H}^s(F)=0\}=\sup\{s\geq 0: \mathcal{H}^s(F)=\infty          \}.
\]

\textbf{Box dimensions:}
The upper/lower box dimension of $F$ is:
\[
\overline{\boxd} \text{ resp. }\underline{\boxd}(F)=\limsup_{r\to 0} \text{ resp. }\liminf_{r\to 0}\left(-\frac{\log N_r(F)}{\log r}\right).
\]
If the limsup and liminf are equal we call this value the box dimension of $F$.

Throughout this paper, we shall discuss cubes or squares. To be precise when we have an Euclidean space $\mathbb{R}^{n}$, we fix a Cartesian coordinate system. A cube centred at $x\in\mathbb{R}^n$ with side length $r$ is the following 'layer' set of the supreme norm in Euclidean space:
\[
\{y\in\mathbb{R}^n: \|y-x\|_{\infty}= r/2\}.
\]
So we see that such a cube is aligned with the coordinate axis.

We shall also consider the Assouad dimension in this paper, see \cite{F} for more details.

\textbf{Assouad dimension:} The Assouad dimension of $F$ is 
\begin{align*}
\Assouad F = \inf \Bigg\{ s \ge 0 \, \, \colon \, (\exists \, C >0)\, (\forall & R>0)\,  (\forall r \in (0,R))\, (\forall x \in F) \\ 
&N_r(B(x,R) \cap F) \le C \left( \frac{R}{r}\right)^s \Bigg\}
\end{align*}
where $B(x,R)$ denotes the closed ball of centre $x$ and radius $R$.

\textbf{Approximation symbols:} When counting covering numbers, it is convenient to introduce notions $\approx, \lessapprox, \gtrapprox$ for approximately equal, approximately smaller and approximately larger.

As our box counting procedure always involves scales, later we use $1>\delta>0$ to denote a particular scale. Then for two quantifies $f(\delta), g(\delta)$ we define the following:
\[
f\lessapprox g\iff \forall \epsilon>0, \exists C_{\epsilon}>0 \text{ such that } \forall \delta>0, f(\delta)\leq C_{\epsilon} \delta^{-\epsilon} g(\delta).
\] 
\[
f\gtrapprox g\iff g\lessapprox f.
\]
\[
f\approx g\iff f\lessapprox g \text{ and } g\lessapprox f.
\]

Thus $f\approx g$ in our box counting procedure can be intuitively read as "$f$ and $g$ give the same box dimension". Later in context, $f,g$ can be either box covering number of scale $\delta$ or Lebesgue measure of a $\delta$ neighbourhood of a set.
\section{Proof of theorem \ref{Th1}}

We shall focus on $\mathbb{R}^2$. in this situation we have a clear picture of what is going on. The arguments work for other cases after some modification. From now on we fix $n=2$ and let $G$ be a set which contains cubes of all side lengths in $(0,1)$.

The analogous result for Hausdorff dimension is obvious since $G$ contains boundary of a cube and therefore has Hausdorff dimension at least $1$.

We now consider the box dimension result. Let $\delta>0$ be a small positive number. Then we can find approximately $\frac{1}{100\delta}$ many $100\delta$-separated points in $[1/2,1]$. We denote those points from small to large as:
\[
r_1<r_2<\dots<r_k, k\approx (100\delta)^{-1}
\] 

For each $r_i,i\in [1,k]$, there is a cube $C_{r_i}$ of side length $r_i$ contained in $G$. Any cube contains $4$ sides, therefore there are $\approx 4k$ many sides of length in $[1/2,2]$. Now we can focus on for example the right most side $I_{r_i}$ of each cube $C_{r_i}$. We know that the obstruction of $G$ having high dimension is the heavy overlap of sides. Then we consider $\delta$-neighbourhood $C^{\delta}_{r_i}$ of each $C_{r_i}$. Consider the characteristic function $\chi_{I^\delta_{r_i}}$. If there exist a $x\in\mathbb{R}^2$ and an integer $M>0$ such that:
\[
\sum_{i} \chi_{I^\delta_{r_i}}(x)\geq M,
\]
then there are $M$ cubes whose right sides are $2\delta$ close to each other. Because different cubes have side length at least $100\delta$ difference, we see that the left sides of those $M$ cubes stays at least $96\delta$ away from each other. This implies that we can find $M$ sides which are $96\delta$ away from each other therefore union of the $\delta$-neighbourhood of those sides takes area at least:
\[
0.5M \delta.
\]

Now let $M$ be such that:
\[
\left\|\sum_{i} \chi_{I^\delta_{r_i}}\right\|_{\infty}= M,
\]

we see that:
\[
\left\| \chi_{\bigcup_i I^\delta_{r_i}}\right\|_{L^1}\geq \frac{1}{M} \left\|\sum_{i} \chi_{I^\delta_{r_i}}\right\|_{L^1}\gtrapprox \frac{0.5 \delta\times \frac{1}{100\delta}}{M}=\frac{1}{200M}.
\]

From the above argument we see that $G^\delta$ takes area at least:
\[
\max\left\{0.5M \delta, \frac{1}{200M}\right\}\gtrapprox \frac{1}{40}\sqrt{\delta}.
\]

This gives us the lower bound of the lower box dimension of $G$:
\[
\underline{\boxd} G\geq 1.5.
\]

For other cases, "sides" of cubes are replaced by "faces" which are subsets of $n-1$ dimensional affine hyperplanes (here we allow the case $n=1$). Arguing as above we see that there is a constant $c>0$:
\[
|G^{\delta}|\geq c\max\left\{M \delta, \frac{1}{M}\right\}\gtrapprox \sqrt{\delta}.
\]
This implies that $\underline{\boxd} G\geq n-0.5$. Therefore we established the lower bound for lower box dimension, the Assouad dimension result follows because the Assouad dimension is always greater or equal to the lower box dimension.

\section{Proof of theorem \ref{Th1.1}}

\textbf{Dimension one:}

 We want to find $F_1,F_2$ such that $\Haus F_1=0$, $\boxd F_2=0.5$ and the distance set $|F_1-F_1|, |F_2-F_2|$ contain $[0,1]$. And for every $\sigma>0$ we want to find $F_3$ such that $\Assouad F_3\leq 1/2+\sigma$ and $|F_3-F_3|$ contains $[0,1]$. 

The construction of $F_1,F_2$ can be found in \cite{D}. For Assouad dimension, we shall use a result in \cite{N} by Nathanson which says that for all integer $n>1$ there exists a subset $B\subset\{0,\dots,n-1\}$ such that $B+B \mod n=\{0,\dots,n-1\}$ and $|B|\leq 2(n\log n)^{1/2}+2$.

Now for any integer $n>1$, we find such set $B$, then we let $C=B\cup (n-B\mod n)$. We can now construct Cantor set in $[0,1]$ by restricting $n$-ary digital expansions. Precisely, we let \[C_n=\{x\in [0,1]: n\text{-ary expansion of $x$ contains only digits in } C \}\]
Then it is easy to see that $|C_n-C_n|$ contains $[0,1]$. We also have the following result concerning the Assouad dimension of $C_n$:
\[
\Assouad C_n=\Haus C_n=\frac{\log |B|}{\log n}\leq \frac{\log (4(n\log n)^{1/2}+4)}{\log n}.
\]
Let $F_3=C_n$ for a large enough $n$ we see that $\Assouad F_3\leq 1/2+\sigma$.

\begin{rem}
	For set $F_3$, we can also apply Solomyak's result \cite{S} on Palis conjecture if we only require $|F_3-F_3|$ to have positive measure.
\end{rem}

\textbf{Higher dimensional cases:}

Having settled the one dimensional case, we shall see that we can extend the one dimensional case to higher dimensional cases. Again we focus here on $\mathbb{R}^2$ and similar arguments lead us the corresponding results in $\mathbb{R}^n,n\geq 3.$

Pick $F\subset [0,1]\times \{0\}\subset\mathbb{R}^2$ whose property we shall require later. Then for $x\in F$ we construct two lines passing through $x$ with slope $\pm 1.$ Denote $L(F)$ be the union of all such lines. Let $r\in |F-F|$, then there exist $x_1,x_2\in F$ with $|x_1-x_2|=r$. Then the lines passing through $x_1,x_2$ with slope $\pm 1$ (there are four of them) will enclose a cube of side length $r/\sqrt{2}$.

Now it is easy to see that $L(F)$ can be written as a union of two subsets with lines of slope $1$ or $-1$. Each of those two sets can be viewed as the Cartesian product of $F$ and $\mathbb{R}$ with a certain affine transformation. From this fact it is easy to see that:
\[
\dim (L(F))=\dim F+1,
\]
here $\dim$ can be any dimension we considered above. For example we see that:
\[
\dim_H L(F_1)=1
\]
and $L(F_1)$ contains cubes of all side length in $[0, \sqrt{2}/2]$. After some rescaling, we can obtain a set of Hausdorff dimension $1$ which contains boundaries of cubes with all side lengths in $(0,1)$.

The box and Assouad dimension results follow in a similar way.

For higher dimensions, instead of drawing lines through points in $F$ we shall draw hyperplanes and this concludes Theorem \ref{Th1.1}.

\section{Largeness in general, proof of theorem \ref{Th2}}

We will work now in $\mathbb{R}^2$ and the argument can be modified to $\mathbb{R}^n$.

Denote $\pi_1,\pi_2$ for the coordinate projections. Let $f_p, G, G_r$ be described as in the statement of this theorem.

We denote the stripes of width $\frac{1}{100}$ as
\[
S_i=[0,1]\times \left(\frac{i}{100},\frac{i+1}{100}\right]
\] 
where $i\in\{0,\dots,99\}$. Later we shall choose thinner and thinner strips.

For each $i$, we want to find cubes that pass through $S_i$. Denote the following set
\[
C_i=\left\{t\in (0,1): \left(\frac{i}{100},\frac{i+1}{100}\right]\subset\left(\pi_2(f_p(t))-\frac{t}{4},\pi_2(f_p(t))+\frac{t}{4}\right)\right\}
\]

For $t\in C_i$, the cube $tK_{1}+f_p(t)$ has two of its sides crossing over the stripe $S_i$. The $t/4$ appeared above ensures that $\frac{t}{2}K_{1}+f_p(t)$ also has two of its sides crossing over the stripe $S_i$. We need this fact to study the set $G_r$ for $r\in (1/2,1)$.

Because $f_p$ is Borel measurable, $C_i,i\in\{0,\dots,99\}$ are all Borel measurable. We also see that because $\cup_{i} C_i=(0,1)$, there exist one $i$ such that $C_i$ has positive Lebesgue measure.

Now we focus on what is happening in stripe $S_i$. Because $C_i$ consists $t$ such that $tK_{1}+f_p(t)$ passes through $S_i$, we see that
\[
\left\{t+\pi_1(f_p(t)): t\in C_i\right\}\times \left(\frac{i}{100},\frac{i+1}{100}\right]\subset G.
\] 

We now consider the scaling ratio $r\in [1/2,1)$ and the set $G_r$. We see that
\[
\left\{rt+\pi_1(f_p(t)): t\in C_i\right\}\times \left(\frac{i}{100},\frac{i+1}{100}\right]\subset G_r.
\]

So we now need to study the set

\[
F_r=\left\{rt+\pi_1(f_p(t)): i\in C_i\right\}.
\]

Consider the following union of line segments

\[
F=\bigcup_{t\in C_t} \{(x,y):x\in (1/2,1), y=tx+\pi_1(f_p(t))\}.
\]

Then $F_r$ is precisely the slicing of $F$ with $\{x=r\}$. By the duality principle \cite[Chapter 12, section 1]{Fa} we see that from the fact that $\Haus C_i=1$
\[
\Haus F_r=1
\]
for Lebesgue almost all $r\in (1/2,1)$. For this argument we need $C_i$ to be a Borel set because Marstrand projection theorem is used in the duality principle.

This implies that

\[
\Haus G_r=2
\]

for Lebesgue almost all $r\in (1/2,1)$. For other values of $r$, we need to consider thinner stripes and stronger passing through conditions. More precisely we replace $1/100$ by $1/k$ with some integer $k\geq 101$, and 
\[
C_i=\left\{t\in (0,1): \left(\frac{i}{k},\frac{i+1}{k}\right]\subset\left(\pi_2(f_p(t))-\frac{t}{\log k},\pi_2(f_p(t))+\frac{t}{\log k}\right)\right\}.
\]

Then we obtain full Lebesgue measure of $r\in \left(\frac{2}{\log k},1\right)$ such that
\[
\Haus G_r=2.
\]

We can choose arbitrarily large $k$ and therefore Theorem \ref{Th2} follows when $n=2$.

For $n\geq 3$, the argument is completely similar. We need to decompose $[0,1]^3$ into thin blocks of width $1/100$
\[
S_{i,j}=[0,1]\times \left(\frac{i}{100},\frac{i+1}{100}\right]\times \left(\frac{j}{100},\frac{j+1}{100}\right]
\] 
from here the rest of the proof should be clear and we omit the full details.

\section{Further comments}

\begin{flushleft}
	\textbf{The centre cube packing}: 
\end{flushleft}

We defined two ways of cube packings in definition \ref{defsi}, definition \ref{defpo}. We only showed the largeness in general for size cube packings. It is natural to consider also the largeness in general for centre cube packings. It turns out that we can follow a similar approach and obtain a similar result. For convenience we put the statement here although we will not give a proof.

\begin{thm*}
	For any $n\geq 2$, consider an measurable centre cube packing $g_p$. Let $G=T(g_p)$ and for each $r\in (0,1)$ we construct the following set
	\[
	G_r=\bigcup_{x\in [0,1]^n} rg_p(x) K_{n-1}+x.
	\]
	Then for Lebesgue almost all $r\in (0,1)$, $\Haus G_r=n$. 
\end{thm*}

\begin{flushleft}
	\textbf{Size cube packings with a subset of sizes}:
\end{flushleft}

It is also natural to consider the size cube packing in the following sense

\begin{defn}
	A restricted size cube packing is defined as a map $f_p$ from a parameter set $E\subset (0,1)$ to $[0,1]^n$:
	\[
	f_p:E\to \mathbb{R}^n
	\]
	
	The image of packing $f_p$ is 
	\[
	T(f_p)=\bigcup_{t\in (0,1)} tK_{n-1}+f_p(t)\subset\mathbb{R}^n.
	\]
	
	We call the packing measurable if $E$ is a Borel set and $f_p$ is Borel measurable.
\end{defn} 

When $E$ is a set of Hausdorff dimension $s\in (0,1)$ similar argument of proof of theorem \ref{Th1}, \ref{Th2} gives us the following result whose proof we omit.

\begin{thm*}
	For any $n\geq 2$, consider an measurable restricted size cube packing $f_p$. Let $G=T(f_p)$ and for each $r\in (0,1)$ we construct the following set
	\[
	G_r=\bigcup_{t\in E} rt K_{n-1}+f_p(t).
	\]
	
	Then $\underline{\boxd} G\geq n-1+\frac{s}{2}$ and for Lebesgue almost all $r\in (0,1)$, $\Haus G_r=n-1+s$. Where $s$ is the Hausdorff dimension of $E$.
\end{thm*}

We shall point out some main steps. The first step is to construct a Frostman measure supported on $E$. Then we can argue as in Theorem \ref{Th2} but with Lebesgue measure replaced by this Frostman measure. In the duality argument, we are in the situation of finding good directions of projecting a set of Hausdorff dimension $s$. The good directions form a large set because the exceptional set has Hausdorff diemnsion at most $s$.

\section{Further questions}
In this section we mention some problems which are related with the topic in this paper and seem to be interesting for further study.

\begin{flushleft}
\textbf{Dimension balancing}
\end{flushleft}

Although we have constructed size cube packing sets in $\mathbb{R}^n$ with Hausdorff dimension $n-1$, box dimension $n-0.5$. Those two values are often not simultaneously attained. As we can see from the sharp example, when the Hausdorff dimension is $n-1$ the upper box dimension is actually $n$.

It is quite possible that the Hausdorff dimension and the upper box dimension for the packing problems considered in this paper have some balancing.

We posed the following question

\begin{ques}
	For a size cube packing $f_p$ and the set $G=T(f_p)$. Is it true that
	\[
	\Haus G+\overline{\boxd} G\geq 2n-1.
	\]
\end{ques}

In \cite{CCHK}, a dense $G_\delta$ set argument was used to show that typically such set $G$ has Hausdorff dimension $n-1$. At the same time we from the same dense $G_\delta$ set argument we can also see that typically such set $G$ has upper box dimension $n$. So the dimension balancing holds at least typically.   

\begin{flushleft}
\textbf{Skeleton packings}
\end{flushleft}

In this paper we only considered the packings of $n-1$ skeletons of $n$ dimensional cubes. It is also interesting to consider $k$-skeletons for $0\leq k\leq n-2$. There are some considerations for Hausdorff dimension of these packing problems using Baire category arguments in \cite{CCHK}, \cite{Th} and \cite{SK}. What about box dimensions? How large they must be? To be precise we can define the following size skeleton packing:

\begin{defn}[Vertex cube packing: different sizes]
	A size cube packing is defined as a map:
	\[
	f_p:(0,1)\to [0,1]^n.
	\]
	
	For an integer $k$ such that $0\leq k\leq n$. The $k$ skeleton image of packing $f_p$ is 
	\[
	T_k(f_p)=\bigcup_{t\in (0,1)} tK_{k}+f_p(t)\subset\mathbb{R}^n.
	\]
	
	We call the packing measurable if $f_p$ is Borel measurable.
\end{defn} 

We can pose the following question:

\begin{ques}
	Let $f_p$ be a size cube packing, then what is the sharp lower bound for
	$
	\underline{\boxd} T_k(f_p)
	$?
\end{ques}
\section{Acknowledgement}
This manuscript was written during the author's stay in Institut Mittag-Leffler. The author thanks Jonathan Fraser, Tam\'{a}s Keleti and Meng Wu for discussions.
\bibliography{squaresone}
\bibliographystyle{amsalpha}

\end{document}